\documentclass[12pt]{iopart}
\usepackage{iopams}
\newtheorem{theorem}{Theorem}[section]
\newtheorem{corollary}[theorem]{Corollary}
\newtheorem{lemma}[theorem]{Lemma}
\newtheorem{proposition}[theorem]{Proposition}

\newtheorem{assumption}[theorem]{Assumption}
\newtheorem{remark}[theorem]{Remark}

\begin{document}
\title[Inverse problems for a reaction-diffusion
system]{Inverse problems for a two by two reaction-diffusion
system using a Carleman estimate with one observation}

\author{Michel Cristofol$^1$,
Patricia Gaitan$^2$ 
 and Hichem Ramoul$^{3,*}$}

\address{$^1$ Laboratoire d'analyse, topologie, probabilit\'{e}s
 CNRS UMR 6632, Marseille, France and Universit\'e Aix-Marseille III}
\address{$^2$ Laboratoire d'analyse, topologie, probabilit\'{e}s
 CNRS UMR 6632, Marseille, France and Universit\'e Aix-Marseille II}
\address{$^3$ Centre universitaire de Khenchela,
 Route de Batna, BP 1252, Libert\'e, 40004 Khenchela, Alg\'erie}
\address {$^*$ This research was supported by the program Tassili 04MDU606}
\ead{\mailto{cristo@cmi.univ-mrs.fr}, \mailto{gaitan@cmi.univ-mrs.fr}, \mailto{RAMOUL2477@yahoo.fr}}

\begin{abstract}
  For a two by two reaction-diffusion system on a bounded domain we give a
simultaneous stability result for one
coefficient and for the initial conditions.
  The key ingredient is a global Carleman-type estimate
  with a single observation acting on a subdomain.
\end{abstract}
\noindent{\it Keywords\/} : Inverse problem, Reaction-diffusion system, Carleman estimate
\ams{35K05, 35K57, 35R30.}
\submitto{Inverse Problems}
\maketitle
\normalsize
\section{Introduction}
This paper is devoted to the simultaneous identification 
of one coefficient and the initial conditions
in a reaction-diffusion system using the least number of 
observations as possible.\\
Let $\Omega \subset \mathbb{R}^n$ be a bounded domain of $\mathbb{R}^n$
with $n \le 3$, (the assumption $n \le 3$ is necessary in order to obtain the appropriate regularity
for the solution using classical Sobolev embedding, see \cite{B:83}).
We denote by $\nu$ the outward unit normal to $\Omega$ on $\Gamma=\partial \Omega$
assumed to be
of class $\mathcal{C}^1$. 
Let $T>0$ and $t_0 \in (0,T)$. We shall use the following notations $Q_0=\Omega \times (0,T)$,
$Q=\Omega \times (t_0,T)$, $\Sigma=\Gamma \times (t_0,T)$  
and $\Sigma_0=\Gamma \times (0,T)$.
We consider the following reaction-diffusion system:
\begin{equation}
\left \{ \begin{array}{lll}
 \label{syst}
  \partial_t u= \Delta u+a(x)u+b(x)v & \mbox{in} & Q_0,\\
  \partial_t v= \Delta v+c(x)u+d(x)v & \mbox{in} & Q_0, \\
   u(t,x)=g(t,x), \; v(t,x)=h(t,x) & \mbox{on} & \Sigma_0 ,\\
   u(0,x)=u_0 \;\; \mbox{ and}\;\; v(0,x)=v_0 & \mbox{in} & \Omega.
\end{array}\right.
\end{equation}
Reaction-diffusion systems are frequently used to model
several physical applications, for example :
in biology and medecine, emergence and growth of cancer and angiogenesis (see \cite{CH:05}), 
in ecology, prey-predator systems, insect dispersal (see \cite{M:89}),
in chemistry, reaction in the presence of diffusion could 
produce spatial pattern of the chemical concentration (see \cite{T:52}).\\ \noindent
Our problem can be stated as follows:\\
Is it possible to determine the coefficient $b(x)$ and the initial conditions
$u_0(x)$, $v_0(x)$ for $x \in \Omega$ from the following measurements:
$$\partial_t v_{|(t_0,T) \times \omega}\ \mbox{ and }
\Delta u(T^{\prime}, \cdot), \ u(T^{\prime}, \cdot), \ v(T^{\prime}, \cdot)
\ \mbox{ in } \Omega \mbox{ for }  T^{\prime}=\frac{t_0+T}{2},$$
where $\omega$ be a subdomain of $\Omega$ ?\\
\noindent
Throughout this paper, let us consider the following set
$$\Lambda(R)=\{\Phi\in L^{\infty }(\Omega );\,\, \|\Phi\|_{L^{\infty }(\Omega )}\leqslant R \},$$
where $R$ is a given positive constant.\\ \noindent
If we assume that ($u_0$, $v_0$) belongs to $(H^2(\Omega))^2$ and $g$, $h$
are sufficiently regular  
(e.g. $\exists \;\epsilon>0$ such that $g, h \in H^1(t_0,T, H^{2+\varepsilon}(\partial \Omega)) \cap
H^2(t_0,T, H^{\varepsilon}(\partial \Omega))$), then (\ref{syst}) admits a
solution in $H^1(t_0,T, H^{2}(\Omega))$ (see \cite{L:68}).
We will later use this regularity result.\par \noindent 
Let $(u,v)$ (resp. ($\widetilde{u}$, $\widetilde{v}$)) be solution of (\ref{syst}) associated to
($a$, $b$, $c$, $d$, $u_0$, $v_0$) 
(resp. ($a$, $\widetilde{b}$, $c$, $d$, $\widetilde{u}_0$, $\widetilde{v}_0$)) 
satisfying some regularity and "positivity" 
properties:
\begin{equation} \label{reg-coeff}
 \left\{ \begin{array}{l}
a,\ b,\ c,\ d,\ \widetilde{b} \in \Lambda(R),\\
\mbox{There exist two constants}\,\, r>0, \,\, c_0>0
\,\,\mbox{such that}\\ \widetilde{u}_0 \geq 0, \,\, \widetilde{v}_0 \geq r,\,\,
 c \geq c_0,\,\, \widetilde{b} > 0, \, \, c + dr \geq 0, \\
 g \geq 0 \,\, \mbox{and} \, \, h \geq r.
\end{array}
\right. 
\end{equation}
Such assumptions allows us to state that the function $\widetilde{v}$
satisfies $| \widetilde{v}(x,T')| \ge r>0$ in $\Omega$ (see \cite{S:83} thm 14.7 p. 200).\\ \noindent
 We assume that we can measure $\partial_t v$ on $\omega$ 
in the time interval $(t_0, T)$ for some $t_0 \in (0, T)$ and
$\Delta u$, $u$ and $v$ in $\Omega$ at time $T' \in (t_0, T)$.
Our main results are
\begin{itemize}
\item A stability result for the coefficient $b(x)$ (or $a(x)$): \\
For $\widetilde{u}_0$, $\widetilde{v}_0$ in $H^2(\Omega)$ there exists a
constant
$$
C=C(\Omega, \omega, c_0, t_0, T, r, R) >0
$$
such that
\begin{eqnarray*}  
    |b-\widetilde{b}|^2_{L^2(\Omega)} \leq
    C |\partial_t v - \partial_t \widetilde{v}|^2_{L^2((t_0,T) \times \omega)}
    + C |\Delta u (T', \cdot)- \Delta \widetilde{u}(T', \cdot)|^2_{L^2(\Omega)}\\ 
    + C | u (T', \cdot)-  \widetilde{u}(T', \cdot)|^2_{L^2(\Omega)}
    + C |v (T', \cdot)- \widetilde{v}(T', \cdot)|^2_{L^2(\Omega)}.
\end{eqnarray*}
\item A stability estimate for the initial conditions $u_0$, $v_0$:\\  
For $u_0$, $v_0$, $\widetilde{u}_0$, $\widetilde{v}_0$ in $H^4(\Omega)$ there
exists a constant
$$
C=C(\Omega, \omega, c_0, t_0, T, r, R)>0
$$
such that
$$|u_0-\widetilde{u}_0|^2_{L^2(\Omega)}+|v_0-\widetilde{v}_0|^2_{L^2(\Omega)}\leq \frac{C}{|\log E|},$$
$$\mbox{where}\,\, E=|\partial_t v - \partial_t
\widetilde{v}|^2_{L^2((t_0,T) \times
\omega)}+|u(T',\cdot)-\widetilde{u}(T',\cdot)|^2_{H^2(\Omega)}
+|v(T',\cdot)-\widetilde{v}(T',\cdot)|^2_{H^2(\Omega)}.$$
\end{itemize}
The key ingredient to these stability results is a global Carleman estimate 
for a two by two system with one observation.
Controllability for such parabolic systems has been studied in 
\cite{ABDK:05}. The Carleman estimate obtained in 
\cite{ABDK:05} cannot be used to solve the inverse problem of identification of
one coefficient and initial conditions
because of the weight functions which are different in the left and right hand
side of their estimate.
We establish a new Carleman estimate with one observation
involving the same weight function in the left and right hand side.
We prove a stability result for the initial conditions following the method of  \cite{YZ:01}.
Concerning the stability of the initial conditions we use 
an extension of the logarithmic convexity method
(see \cite{I:98}).
We can also cite \cite{K:06}; he provides an estimate for the
initial condition for a general parabolic operator, while the logarithmic convexity method
works only for self-adjoint operators (see section 2 and Theorem 3 in \cite{K:06}).\\
The simultaneous reconstruction of one coefficient
and initial conditions from the measurement of one solution $v$
over $(t_0,T) \times \omega$ and some measurement at fixed
time $T'$ is an essential aspect of our result.
In the perspective of numerical reconstruction,
such problems are ill-posed.
Stability results are thus of importance.\\
For the first time, the method of Carleman estimates was introduced in the field of inverse problems
in the work of Bukhgeim and Klibanov \cite{BK:81}; also see, e.g., \cite{B:99}, \cite{K:84} and \cite{K:92}
 for some follow up publications of
these authors. So far, the method of \cite{BK:81} is the only one enabling to prove uniqueness and stability results
for inverse problems with single measurement data in the $n$-dimensional case with $n \ge 2$, which has generated many
publications, including this one. While this method can provide H\"{o}lder stability results, the topic of the
Lipschitz stability is a more delicate one. The first Lipschitz stability result for a multidimensional inverse
problem (for a hyperbolic equation) was obtained by Puel and Yamamoto \cite{PY:96}, using a modification of the
idea of \cite{BK:81}. In \cite{IY:98} this result was extended for a parabolic equation. 
The main difference between our work
and \cite{IY:98} is that we consider a coupled system of parabolic equations, and the additional data are given only
for one component of this system, the function $\partial_t v(x, t)$ for $(x, t) \in \omega \times (t_0, T)$ .
Inverse problems for parabolic equations
are well studied (see \cite{CD:98}, \cite{IY:98}, \cite{YZ:01}).
A recent book of Klibanov and Timonov \cite{KT:04} is devoted to
the Carleman estimates applied to inverse coefficient problems.
In our knowledge, there is no work about inverse problems 
for coupled parabolic systems.\\
The used method allows us to give a stability result for the
coefficient $a(x)$ adapting assumption \ref{ass1}.
On the other hand, since we only measure $\partial_t v$ on $\omega$,
 we cannot obtain such stability results
for the coefficients $c(x)$ or $d(x)$ of the second equation of (\ref{syst}).
For the reconstruction of two coefficients the problem is more
complicated.
We obtain partial results with restrictive assumptions
on the coefficients $a(x)$, $b(x)$, $c(x)$ and $d(x)$.
In order to avoid such assumptions, we think
it is necessary to use other methods such as those used in \cite{IIY:03}.\\
Our paper is organized as follows. 
In Section 2, we derive a global Carleman estimate 
for system (\ref{syst}) 
with one observation, i.e. the measurement of one solution $v$
over $(t_0,T) \times \omega$.
In Section 3, we prove a stability result for the coefficient $b(x)$ 
when one of the solutions 
$\widetilde{v}$ is in a particular class of solutions with some 
regularity and "positivity" properties.
In Section 4, we prove a stability result for the initial conditions.
%
%
\section{Carleman estimate}
\label{sec:1}
We prove here a Carleman-type estimate with a single observation 
acting on a subdomain $\omega$ of $\Omega$ 
in the right-hand side of the estimate.  Let us
introduce the following notations:
 let $\omega' \Subset \omega$ and let $\widetilde{\beta}$ be a 
 ${\cal C}^2(\Omega)$ function such that
$$
\widetilde{\beta}>0,\,\,\mbox{in}\,\,\Omega,\,\,\widetilde{\beta}=0\,\,\mbox{on}\,\,\partial
\Omega,\,\, \min\{|\nabla \widetilde{\beta}(x)|,\; x \in \overline{\Omega
\setminus \omega'}\}>0 \;\; \mbox{ and }\,\,{\partial}_{\nu}
\widetilde{\beta}< 0\;\;\mbox{on}\;\;\partial \Omega.
$$
Then, we define $\beta= \widetilde{\beta}+K$ with
$K= m \|\widetilde{\beta}\|_{\infty}$ and $m>1$. For $\lambda> 0$ and $t \in
(t_0,T)$, we define the following weight functions
$$
  \label{wf}
  \varphi(x,t)=\frac{e^{\lambda \beta(x)}}{(t-t_0)(T-t)},
  \quad \quad \eta(x,t)=\frac{e^{2\lambda K} -e^{\lambda
  \beta(x)}}{(t-t_0)(T-t)}.
$$
If we set $\psi=e^{-s \eta}q$, we also introduce the following
operators
\begin{eqnarray*}
  M_1\psi  & =& -\Delta\psi - s^{2}\lambda^{2}|\nabla\beta|^2\varphi^{2}\psi + s(\partial_t{\eta})\psi, \\
  M_2\psi  & =& \partial_t \psi +2 s \lambda \varphi \nabla \beta . \nabla \psi +2
s \lambda^{2} \varphi | \nabla \beta |^2 \psi.
\end{eqnarray*}
Then the following result holds (see \cite{FG:05}).\\
\begin{theorem}
\label{th-Carl-Fur} There exist $\lambda_0=\lambda_0(\Omega,
\omega)\geq 1$, $s_0=s_0(\lambda_0, T)>1$ and a positive
constant $C_0=C_0(\Omega, \omega, T)$ such that, for any $\lambda \ge
\lambda_0$ and any $s \ge s_0 $, the following inequality holds:
\begin{eqnarray}
\label{Carl-Fur1}
    \|M_1(e^{-s \eta}q)\|^2_{L^2(Q)} + \|M_2(e^{-s \eta}q)\|^2_{L^2(Q)}\\
  +s \lambda^2 \int \hspace{-6.5pt} \int_Q e^{-2s \eta} \varphi |\nabla q|^2 \ d x\ d t
  +s^{3} \lambda^{4}\int \hspace{-6.5pt} \int_Q e^{-2s \eta} \varphi^{3} |q|^2\ d x\ d t  \nonumber\\
  \leq C_0  \left[
   s^{3} \lambda^{4}\int_{t_0}^T \hspace{-6.5pt}
\int_\omega e^{-2s \eta} \varphi^{3} |q|^2\ d x\ d t
  +\int \hspace{-6.5pt} \int_Q e^{-2s \eta}\ | \partial_t q - \Delta q|^2\ d x\ d t\right],\nonumber
\end{eqnarray}
for all $q \in H^1(t_0,T, H^2(\overline{\Omega}))$ with $q=0$ on
$\Sigma$.
\end{theorem}
From the above theorem we have also the following result (see
\cite{F:00} and \cite{IY:98}).\\
\begin{proposition}
\label{pr-Carl-Fur} There exist $\lambda_0=\lambda_0(\Omega,
\omega)\geq 1$, $s_0=s_0(\lambda_0, T)>1$ and a positive
constant $C_0=C_0(\Omega, \omega, T)$ such that, for any $\lambda \ge
\lambda_0$ and any $s \ge s_0 $, the following inequality holds:
\begin{eqnarray}
\label{Carl-Fur2}
    s^{-1} \int \hspace{-6.5pt} \int_Q e^{-2s \eta} \varphi^{-1}
(|\partial_t q|^2+|\Delta q|^2) \ d x\ d t \\
  +s \lambda^2 \int \hspace{-6.5pt} \int_Q e^{-2s \eta} \varphi |\nabla q|^2 \ d x\ d t
  +s^{3} \lambda^{4}\int \hspace{-6.5pt} \int_Q e^{-2s \eta} \varphi^{3} |q|^2\ d x\ d t \nonumber \\
  \leq C_0  \left[
   s^{3} \lambda^{4}\int_{t_0}^T \hspace{-6.5pt}
\int_\omega e^{-2s \eta} \varphi^{3} |q|^2\ d x\ d t
  +\int \hspace{-6.5pt} \int_Q e^{-2s \eta}\ | \partial_t q - \Delta q|^2\ d x\ d t\right],\nonumber
\end{eqnarray}
for all $q \in H^1(t_0,T, H^2(\overline{\Omega}))$ with $q=0$ on
$\Sigma$.
\end{proposition}
We consider the solutions ($u$, $v$) and ($\widetilde{u}$, $\widetilde{v}$)
to the following systems
\begin{equation}
\left \{ \begin{array}{lll}
  \label{syst1}
    {\partial}_t u= \Delta u+au+bv & \mbox{in} & Q_0,\\
    {\partial}_t v= \Delta v+cu+dv & \mbox{in} & Q_0, \\
    u(t,x)=g(t,x), \; v(t,x)=h(t,x) & \mbox{on} & \Sigma_0,\\
    u(0,x)=u_0 \;\; \mbox{and} \;\; v(0,x)=v_0 & \mbox{in} & \Omega, 
\end{array}\right.
\end{equation}
and
\begin{equation}
\left \{ \begin{array}{lll}
  \label{syst2}
    {\partial}_t \widetilde{u}= \Delta \widetilde{u}+a\widetilde{u}+\widetilde{b}\widetilde{v} & \mbox{in} & Q_0,\\
    {\partial}_t \widetilde{v}= \Delta \widetilde{v}+c\widetilde{u}+d\widetilde{v} & \mbox{in} & Q_0,\\
    \widetilde{u}(t,x)=g(t,x), \; \widetilde{v}(t,x)=h(t,x) & \mbox{on} & \Sigma_0,\\
    \widetilde{u}(0,x)=\widetilde{u}_0 \;\; \mbox{and} \;\; \widetilde{v}(0,x)=\widetilde{v}_0 & \mbox{in} & \Omega.
\end{array}\right.
\end{equation}
We set $U=u-\widetilde{u}$, $V=v-\widetilde{v}$, $y=\partial_t(u-\widetilde{u})$, 
$z=\partial_t(v-\widetilde{v})$ and $\gamma=b-\widetilde{b}$.
Then $(y,z)$ is solution to the following problem
\begin{equation}
\left \{ \begin{array}{lll}
  \label{syst3}
    {\partial}_t y= \Delta y+ay+bz+\gamma \partial_t \widetilde{v} & \mbox{in} & Q_0,\\
    {\partial}_t z= \Delta z+cy+dz & \mbox{in} & Q_0, \\
    y(t,x)=z(t,x)=0  & \mbox{on} & \Sigma_0,\\
    y(0,x)=\Delta U(0,x)+a U(0,x)+b V(0,x)+ \gamma \widetilde{v}(0,x), & \mbox{in} & \Omega,\\
    z(0,x)=\Delta V(0,x)+cU(0,x)+dV(0,x) & \mbox{in} & \Omega.
\end{array}\right.
\end{equation}
Note that the previous initial conditions are available for all $T' \in (0,T)$.
Indeed with (\ref{syst1}) and (\ref{syst2}) we can determine $y(T',x)$, $z(T',x)$ and we obtain
$$y(T',x)=\Delta U(T',x)+aU(T',x)+bV(T',x)+\gamma \widetilde{v}(T',x),$$
$$z(T',x)=\Delta V(T',x)+c U(T',x)+d V(T',x).$$
We consider the functional
\begin{eqnarray*}
\label{funct}
    I(q)= s^{-1} \int \hspace{-6.5pt} 
\int_Q e^{-2s \eta} \varphi^{-1} (|\partial_t q|^2+|\Delta q|^2) 
  \ d x\ d t \\
  + s \lambda^2 \int \hspace{-6.5pt} \int_Q e^{-2s \eta} \varphi |\nabla q|^2 \ d x\ d t 
  +s^{3} \lambda^{4}\int \hspace{-6.5pt} \int_Q e^{-2s \eta} \varphi^{3} |q|^2\ d x\ d t.\nonumber 
\end{eqnarray*}
Then using the Carleman estimate (\ref{Carl-Fur2}), the solution ($y$, $z$) of (\ref{syst3})
satisfies
\begin{eqnarray}
\label{est1}
   I(y)+I(z) \leq C_1 [
   s^{3} \lambda^{4}\int_{t_0}^T \hspace{-6.5pt} 
\int_\omega e^{-2s \eta} \varphi^{3} |z|^2\ d x\ d t \\
   + s^{3} \lambda^{4}\int_{t_0}^T \hspace{-6.5pt} 
\int_\omega e^{-2s \eta} \varphi^{3} |y|^2\ d x\ d t \nonumber \\ 
   +\int \hspace{-6.5pt} \int_Q e^{-2s \eta}\ 
(| ay|^2+|bz|^2 + |\gamma \partial_t \widetilde{v}|^2)\ d x\ d t
  +\int \hspace{-6.5pt} \int_Q e^{-2s \eta}\ (| cy|^2+|dz|^2)\ d x\ d t ]\nonumber 
\end{eqnarray}
Let $\xi$ be a smooth cut-off function satisfying
$$
\left \{ \begin{array}{ll}
\xi(x)=1  & \forall x \in \omega',\\
0<\xi(x) \leq 1 & \forall x \in \omega'',\\
\xi(x)=0 & \forall x \in \mathbb{R}^n \setminus \omega'',
\end{array}\right.
$$
\noindent
where $\omega' \Subset \omega'' \Subset \omega \Subset \Omega $.\\
\noindent
We shall estimate the following three terms 
$$I:=  s^{3} \lambda^{4}\int_{t_0}^T \hspace{-6.5pt} \int_\omega  \, 
e^{-2s \eta} \varphi^{3} |y|^2\ d x\ d t,$$
$$J:= \int \hspace{-6.5pt} \int_Q e^{-2s \eta}\ | bz|^2\ d x\ d t\ \ 
\mbox {or} \ \ \int \hspace{-6.5pt} \int_Q e^{-2s \eta}\ | dz|^2\ d x\ d t,$$
$$K:= \int \hspace{-6.5pt} \int_Q e^{-2s \eta}\ | ay|^2\ d x\ d t \ \ 
\mbox {or}\ \ \int \hspace{-6.5pt} \int_Q e^{-2s \eta}\ | cy|^2\ d x\ d t.$$
The main difficulty is in estimating from the above the term $I$, because values of the function $y(x,t)$
are unknown for $(x,t) \in \omega \times(t_0,T)$. In doing so, we rely on the assumption $c(x) \ge c_0$ where
$c_0$ is a positive constant, see (\ref{est-I}).
For the first term $I$, we multiply the second equation of (\ref{syst3}) by $s^3 \lambda^4 e^{-2s \eta} \xi \varphi^3 y$
and we integrate over $(t_0,T) \times \omega$. We obtain
$$I':= s^{3} \lambda^{4} \int_{t_0}^T \hspace{-6.5pt} 
\int_\omega c \, e^{-2s \eta} \xi \varphi^{3} |y|^2\ d x\ d t =
  s^{3} \lambda^{4}\int_{t_0}^T \hspace{-6.5pt} 
\int_\omega e^{-2s \eta} \xi \varphi^{3} ( {\partial}_t z - \Delta z -d z) y \ d x\ d t$$
 $$ = s^{3} \lambda^{4}\int_{t_0}^T \hspace{-6.5pt} 
\int_\omega e^{-2s \eta}\xi \varphi^{3} ({\partial}_t z)  y \ d x\ d t
  - s^{3} \lambda^{4}\int_{t_0}^T \hspace{-6.5pt} \int_\omega e^{-2s \eta} \xi
  \varphi^{3} (\Delta z)  y \ d x\ d t$$
  $$- s^{3} \lambda^{4}\int_{t_0}^T \hspace{-6.5pt} 
\int_\omega d e^{-2s \eta} \xi \varphi^{3}  z  y \ d x\ d t=I_1+I_2+I_3.$$
  By integration by parts with respect to the time variable, 
the first integral, $I_1$, can be written as
 \begin{eqnarray*} 
I_1  =  -s^{3} \lambda^{4}\int_{t_0}^T \hspace{-6.5pt} 
\int_\omega e^{-2s \eta} \xi \varphi^{3} z (\partial_t y)\ d x\ d t 
  + 2s^{4} \lambda^{4}\int_{t_0}^T \hspace{-6.5pt} \int_\omega e^{-2s \eta} 
  \xi \varphi^{3} (\partial_t \eta)\ z y\ d x\ d t\\
  -3s^{3} \lambda^{4}\int_{t_0}^T \hspace{-6.5pt} \int_\omega   
  e^{-2s \eta} \xi \varphi^{2} (\partial_t \varphi)\,  z y\ d x\ d t.
  \end{eqnarray*}
  We write $ I_1 =  I_1^1+I_1^2$ with
$$I_1^1=-s^{3} \lambda^{4}\int_{t_0}^T \hspace{-6.5pt} 
\int_\omega e^{-2s \eta} \xi \varphi^{3} z (\partial_t y)\ d x\ d t,$$
$$I_1^2= 2s^{4} \lambda^{4}\int_{t_0}^T \hspace{-6.5pt} \int_\omega e^{-2s \eta} 
  \xi \varphi^{3} (\partial_t \eta)\ z y\ d x\ d t\\
  -3s^{3} \lambda^{4}\int_{t_0}^T \hspace{-6.5pt} \int_\omega   
  e^{-2s \eta} \xi \varphi^{2} (\partial_t \varphi)\,  z y\ d x\ dt .$$
  Using Young inequality, we estimate the two integrals $I_1^1$ and $I_1^2$.
We have
 \begin{eqnarray*}
  |I_1^1 | \leq  s^{3} \lambda^{4} 
  \left[
  C_\varepsilon s^{4} \lambda^{4}\int_{t_0}^T \hspace{-6.5pt} 
\int_\omega e^{-2s \eta} \xi^2 \varphi^{7} |z|^2\ d x\ d t
+ \varepsilon s^{-4} \lambda^{-4}\int_{t_0}^T \hspace{-6.5pt} 
\int_\omega  e^{-2s \eta} \varphi^{-1} 
|\partial_t y|^2\ d x\ d t \right]\\
\leq  
\frac{1}{2 \varepsilon} s^{7} \lambda^{8}\int_{t_0}^T \hspace{-6.5pt} \int_\omega e^{-2s \eta} 
\varphi^{7} |z|^2\ d x\ d t
+ \frac{\varepsilon}{2} s^{-1}  \int \hspace{-6.5pt} \int_Q  e^{-2s \eta} \varphi^{-1} 
|\partial_t y|^2 \ d x \ d t. 
\end{eqnarray*}
The last term of the previous inequality can be "absorbed"
 by the terms in $I(y)$ for $\varepsilon$ sufficiently small
(for exemple $\varepsilon = \frac{1}{2}$).
\begin{eqnarray*}
  |I_1^2 | \leq   C s^{4} \lambda^{4} \left[
  s \lambda \int_{t_0}^T \hspace{-6.5pt} 
\int_\omega  e^{-2s \eta} \xi^2 \varphi (\varphi^2 |\partial_t \eta|^2 
+ |\partial_t \varphi|^2) |z|^2\ d x \ d t \right.\\
  \left. + s^{-1} \lambda^{-1}\int_{t_0}^T \hspace{-6.5pt} 
\int_\omega  e^{-2s \eta} \varphi^{3} |y|^2\ d x\ d t \right]\\
\leq C \left [s^{5} \lambda^{5}\int_{t_0}^T \hspace{-6.5pt} 
\int_\omega e^{-2s \eta} \varphi^{7} |z|^2\ d x\ d t
+ s^{3} \lambda^{3}\int \hspace{-6.5pt} 
\int_Q  e^{-2s \eta} \varphi^{3} |y|^2\ d x\ d t \right].
\end{eqnarray*}
The last inequality holds through the following estimates 
$$ |\partial_t \varphi| \leq C(\Omega, \omega) T \varphi^2, \quad |\partial_t \eta| \leq C(\Omega, \omega) T
\varphi^2, \quad \varphi \leq C(\Omega, \omega) T^4 \varphi^3. $$
The last term of the previous inequality can be "absorbed" by 
the terms in $I(y)$ for $s$ and $\lambda$ sufficiently large.
Finally, we obtain
$$|I_1| \leq C s^{7} \lambda^{8}\int_{t_0}^T \hspace{-6.5pt} 
\int_\omega e^{-2s \eta} \varphi^{7} |z|^2\ d x\ d t+ \mbox{"absorbed terms" },$$
where $C$ is a generic constant which depends on $\Omega$, $\omega$ and $T$.\\
Integrating by parts the second integral $I_2$ 
with respect to the space variable, we obtain
$$I_2= -s^{3} \lambda^{4}\int_{t_0}^T \hspace{-6.5pt} 
\int_\omega \Delta ( e^{-2s \eta} \xi \varphi^{3} y ) z \ d x\ d t.$$
If we denote by $P=  e^{-2s \eta} \xi \varphi^{3}$, then we have
$$I_2=-s^{3} \lambda^{4}\int_{t_0}^T \hspace{-6.5pt} \int_\omega (P \Delta y +2 \nabla P \nabla y + y \Delta P) z \ d x\ d t.$$
We compute $\nabla P$ and $\Delta P$ and we obtain the following estimation for $I_2$
\begin{eqnarray*}
|I_2| \leq  s^{3} \lambda^{4} \biggl[ \varepsilon s^{-4} \lambda^{-4} 
\int_{t_0}^T \hspace{-6.5pt} \int_\omega  e^{-2s \eta} \varphi^{-1} |\Delta y|^2\ d x\ d t 
+ C_\varepsilon s^{4} \lambda^{4} \int_{t_0}^T \hspace{-6.5pt} 
\int_\omega  e^{-2s \eta} \varphi^{7} |z|^2\ d x\ d t \\
+ \varepsilon s^{-2} \lambda^{-2}\int_{t_0}^T \hspace{-6.5pt} 
\int_\omega  e^{-2s \eta} \varphi |\nabla y|^2\ d x\ d t 
+C_\varepsilon s^{2} \lambda^{2} \int_{t_0}^T \hspace{-6.5pt} 
\int_\omega  e^{-2s \eta} \varphi^{5} |z|^2\ d x\ d t \\
 + \varepsilon \int_{t_0}^T \hspace{-6.5pt} 
\int_\omega  e^{-2s \eta} \varphi^{3} |y|^2\ d x\ d t 
+ C_\varepsilon \int_{t_0}^T \hspace{-6.5pt} 
\int_\omega  e^{-2s \eta} \varphi^{3} |z|^2\ d x\ d t \biggr].
\end{eqnarray*}
Therefore we obtain
\begin{eqnarray*}
|I_2| \leq  \varepsilon \biggl[  s^{-1} 
\int \hspace{-6.5pt} \int_Q  e^{-2s \eta} \varphi^{-1} |\Delta y|^2\ d x\ d t
+ s \lambda^{2} \int \hspace{-6.5pt} 
\int_Q  e^{-2s \eta} \varphi |\nabla y|^2\ d x\ d t\\
+ s^{3} \lambda^{4}\int \hspace{-6.5pt} 
\int_Q  e^{-2s \eta} \varphi^{3} |y|^2\ d x\ d t \biggr] 
+ C_\varepsilon \biggl[ s^{7} \lambda^{8}
\int_{t_0}^T \hspace{-6.5pt} \int_\omega  e^{-2s \eta} \varphi^{7} |z|^2\ d x\ d t \\
+s^{5} \lambda^{6}\int_{t_0}^T \hspace{-6.5pt} 
\int_\omega  e^{-2s \eta} \varphi^{5} |z|^2\ d x\ d t
+ s^{3} \lambda^{4} \int_{t_0}^T \hspace{-6.5pt} 
\int_\omega  e^{-2s \eta} \varphi^{3} |z|^2\ d x\ d t \biggr].
\end{eqnarray*}
The first three integrals of the r.h.s. of the previous inequality can be "absorbed" 
by the terms in $I(y)$ for $\varepsilon$ sufficiently small.
Finally, we have
$$|I_2| \leq  C s^{7} \lambda^{8} \int_{t_0}^T \hspace{-6.5pt} 
\int_\omega  e^{-2s \eta} \varphi^{7} |z|^2\ d x\ d t + \mbox{ "absorbed terms" }.$$
For the last integral $I_3$, we have
\begin{eqnarray*}
  |I_3 |
  \leq  C s^{3} \lambda^{4} \left[
   C_\varepsilon \int_{t_0}^T \hspace{-6.5pt} 
\int_\omega e^{-2s \eta} \varphi^{3} |z|^2\ d x\ d t
+ \varepsilon \int \hspace{-6.5pt} \int_Q  e^{-2s \eta} \varphi^{3} |y|^2\ d x\ d t \right]
\end{eqnarray*}
Finally, if we assume that there exists a constant $c_0>0$ such that
$c\geq c_0$ in $\omega$, we have thus obtained for $\lambda$ and $s$ sufficiently large and
$\varepsilon$ sufficiently small the following estimate:
\begin{equation} \label{est-I}
|I| \leq \frac{1}{c_0}|I'|\leq \frac{C}{c_0} s^{7} \lambda^{8}
\int_{t_0}^T \hspace{-6.5pt} \int_\omega e^{-2s \eta} \varphi^{7}
|z|^2\ d x\ d t .
\end{equation}
For the integrals $J$ and $K$, since $a,\ b,\ c,\ d\ \in \Lambda(R)$ and using the estimate
$$1\leq C(\Omega, \omega)T^6\varphi^3 /4,$$
we have
$$|J| \leq C \int \hspace{-6.5pt} \int_Q e^{-2s \eta}\ \varphi^3\ |z|^2\ d x\ d t,$$
$$|K| \leq C \int \hspace{-6.5pt} \int_Q e^{-2s \eta}\ \varphi^3\ |y|^2\ d x\ d t,$$
and these terms can be "absorbed" by the terms $I(y)$ and $I(z)$ for $\lambda$
and $s$ sufficiently large. 
If we now come back to inequality (\ref{est1}), using the estimates for $I$,
$J$ and $K$, and choosing $\lambda$
and $s$ sufficiently large and $\varepsilon$ sufficiently small, we can thus write
\begin{eqnarray*}
   I(y)+I(z) \leq C_1 [
   s^{3} \lambda^{4}\int_{t_0}^T \hspace{-6.5pt} 
\int_\omega e^{-2s \eta} \varphi^{3} |z|^2\ d x\ d t 
   + s^{7} \lambda^{8}\int_{t_0}^T \hspace{-6.5pt} 
\int_\omega e^{-2s \eta} \varphi^{7} |z|^2\ d x\ d t \\ 
   +\int \hspace{-6.5pt} \int_Q e^{-2s \eta}\ | \gamma \partial_t \widetilde{v} |^2\ d x\ d t].
\end{eqnarray*}
Observing that
 $$\varphi^3 \leq C(\Omega, \omega) T^8 \varphi^7,$$
We have thus obtained the fundamental result
\begin{theorem}
\label{th-carl-est} We assume $a,\ b,\ c,\ d\ \in \Lambda(R)$ and that
exists $c_0>0$ such that $c\geq c_0$ in $\omega$. Then there
exist $\lambda_1=\lambda_1(\Omega, \omega)\geq 1$,
$s_1=s_1(\lambda_1, T)>1$ and a positive constant $C_1=C_1(\Omega,
\omega, c_0, R, T)$ such that, for any $\lambda \geq \lambda_1$ and
any $s \geq s_1 $, the following inequality holds:
\begin{equation}
\label{carl-est}
I(y)+I(z) \leq C_1 \left[
    s^{7} \lambda^{8}\int_{t_0}^T \hspace{-6.5pt} \int_\omega e^{-2s \eta} \varphi^{7} |z|^2\ d x\ d t  
   +\int \hspace{-6.5pt} \int_Q e^{-2s \eta}\ | \gamma \partial_t \widetilde{v} |^2\ d x\ d t \right],
   \end{equation}
    for any solution ($y$, $z$) of (\ref{syst3}).
   \end{theorem}
%
%
\section{ Uniqueness and stability estimate with one observation}
%
%
In this section, we establish, a stability
inequality and deduce a uniqueness result for the
coefficient $b$.  This inequality (\ref{res:sta}) estimates the difference between the
 coefficients $b$ and $\widetilde{b}$  with an upper bound given by some Sobolev
norms of the difference between the solutions $v$,  and $\widetilde{v}$ of (\ref{syst1}) and
(\ref{syst2}). 
Recall that $U=u-\widetilde{u}$, $V=v-\widetilde{v}$, $y=\partial_t(u-\widetilde{u})$, 
$z=\partial_t(v-\widetilde{v})$, $\gamma=b-\widetilde{b}$ and
$$
\left \{ \begin{array}{lll}
  \label{syst4}
    {\partial}_t y= \Delta y+ay+bz+\gamma \partial_t \widetilde{v} & \mbox{in} & Q_0,\\
    {\partial}_t z= \Delta z+cy+dz & \mbox{in} & Q_0,\\
    y(t,x)=z(t,x)=0  & \mbox{on} & \Sigma_0,\\
    y(0,x)=\Delta U(0,x)+a U(0,x)+b V(0,x)+ \gamma \widetilde{v}(0,x), 
& \mbox{in} & \Omega,\\
    z(0,x)=\Delta V(0,x)+cU(0,x)+dV(0,x) & \mbox{in} & \Omega.
\end{array}\right.
$$
The Carleman estimate (\ref{carl-est}) proved in the previous section will be the key
ingredient in the proof of such a  stability estimate.

Let $T'= \frac{1}{2}(T+t_0)$ the point for which
$\Phi(t)=\frac{1}{(t-t_0)(T-t)}$ has its minimum value. 
For ($\widetilde{u}$, $\widetilde{v}$) solutions of (\ref{syst2}), we make the following assumption:
\begin{assumption}
\label{ass1}
There exist $r>0$,\,\,$c_0>0$ such that $\widetilde{b} \ge 0,\,\, c \ge c_0 ,\,\, c+dr \ge 0,\,\, \widetilde{u}_0 \ge 0,\,\, \widetilde{v}_0 \ge
r,\,\, g \ge 0 $ and $ h \ge r$.
\end{assumption}
Such assumption allows us to state that the solution $\widetilde{v}$ is such that 
$|\widetilde{v}(x,T')| \geq r > 0$ in $\Omega$ (see \cite{S:83}, theorem 14.7 p.200).
Furthermore if we assume that $\widetilde{u}_0$,
$\widetilde{v}_0$ in $H^2(\Omega)$,the solutions of (\ref{syst2}) belong to $H^1(t_0,T, H^2(\Omega))$.
Then using classical Sobolev embedding (see \cite{B:83}), we can write for $n \le 3$,
that $\partial_t \widetilde{v}$ belongs to $L^2(t_0,T,L^\infty(\Omega))$
and we assume that $|\partial_t \widetilde{v}|_{L^2(t_0,T)} \in \Lambda(R)$.\\
We set $\psi = e^{-s \eta} y$. 
With the operator
\begin{equation}
\label{eq: M2} 
M_2 \psi= \partial_t \psi +2 s \lambda \varphi \nabla \beta . \nabla \psi
+2 s \lambda^2 \varphi | \nabla \beta |^2 \psi,
\end{equation}
we introduce, following \cite{BP:02}, 

\begin{eqnarray*}
  \mathcal{I} = \Re \int_{t_0}^{T'}\hspace*{-9pt}
  \int_{\Omega} M_2 \psi\;\psi\;dx dt
\end{eqnarray*}
We have the following estimates.
\begin{lemma}
  \label{lemma1}
  Let $\lambda \geq \lambda_1$ and $s\geq s_1$ and let $a,\ b,\ c,\ d\ \in \Lambda(R)$.
  We assume that assumption \ref{ass1} is satisfied then there exists a constant $C=C(\Omega, \omega, T)$ such that
  \begin{eqnarray*} 
    |\mathcal{I}| \leq C s^{-3/2} \lambda^{-2}
    \left[s^7 \lambda^{8}  \int_{t_0}^{T}\hspace*{-9pt}
  \int_{\omega} e^{-2s \eta} \varphi^{7}
    | z |^2\ d x\ d t
      + \int_{t_0}^{T}\hspace*{-9pt}
  \int_{\Omega} 
    e^{-2s \eta} |\gamma|^2 | \partial_t \widetilde{v}|^2 d x\ d t\right].
  \end{eqnarray*}
\end{lemma}
{\bf Proof}:\\
  Observe that
  \begin{eqnarray*}
  |\mathcal{I}| \leq   s^{-3/2} \lambda^{-2}
 \left( \int_{t_0}^{T'}\hspace*{-9pt} \int_{\Omega}
  |M_2 \psi|^2 \ d x\ d t \right)^{1/2}
 \left( s^{3} \lambda^{4} \int_{t_0}^{T'}\hspace*{-9pt} \int_{\Omega}
 e^{-2s \eta}  |y|^2 \ d x\ d t \right)^{1/2},
 \end{eqnarray*}
 thus using Young inequality and the estimate $1\leq C^{\prime}T^6\varphi^3$, we obtain
  \begin{eqnarray*}
    |\mathcal{I}| \leq   C s^{-3/2} \lambda^{-2}
    \left(|M_2 \psi|_{L^2(Q)}^2 +
    s^3 \lambda^{4} \int_{t_0}^{T}\hspace*{-9pt} \int_{\Omega}
     e^{-2s\eta} \varphi^3 |y|^2  dx dt\right),
  \end{eqnarray*}
  which yields the result from Carleman estimate (\ref{carl-est}).
\begin{flushright}
\rule{.05in}{.05in}
\end{flushright}
\begin{lemma} 
  \label{lemma2}
  Let $\lambda \geq \lambda_1$, $s\geq s_1$ and let $a,\ b,\ c,\ d\ \in \Lambda(R)$. 
Furthermore,
we assume that $\widetilde{u}_0$, $\widetilde{v}_0$ in $H^2(\Omega)$ and the
assumption \ref{ass1} is satisfied. Then there exists a constant $C=C(\Omega, \omega, T)$ such that
  \begin{eqnarray}
  \label{eq:lemma2}
    \int_{\Omega} 
    e^{{-2s \eta}(T',x)}\
    |\gamma \widetilde{v}(T',x)|^2 \; dx\\
    \leq C s^{-3/2} \lambda^{-2} \left[
   s^7 \lambda^{8}  \int_{t_0}^{T}\hspace*{-9pt}
  \int_{\omega} e^{-2s \eta} \varphi^{7}
    | z |^2\ d x\ d t
     + \int_{t_0}^{T}\hspace*{-9pt}
  \int_{\Omega} 
    e^{-2s \eta} |\gamma|^2 | \partial_t \widetilde{v}|^2 d x\ d t\right]\nonumber \\
    + C  \int_{\Omega} 
    e^{{-2s \eta}(T',x)}\
    |\Delta\ U(T',x) +a U(T',x) + b V(T',x)|^2 \; dx.\nonumber 
\end{eqnarray}
\end{lemma}
{\bf Proof}:\\
We evaluate integral $\mathcal{I}$ using (\ref{eq: M2})
\begin{eqnarray*}
  \mathcal{I}=
  \frac{1}{2} \int_{t_0}^{T'}\hspace*{-9pt}\int_{\Omega} \partial_t |\psi|^2 \; dx dt
  + s \lambda \int_{t_0}^{T'}\hspace*{-9pt}\int_{\Omega}
  {\varphi} \nabla \beta \cdot \nabla |\psi|^2 \; dx dt
  + 2s \lambda^2 \int_{t_0}^{T'}\hspace*{-9pt}\int_{\Omega} {\varphi}
  \vert \nabla \beta \vert^2 |\psi|^2 dx dt \\
  =\frac{1}{2} \int_{t_0}^{T'}\hspace*{-9pt}\int_{\Omega} \partial_t |\psi|^2 \; dx dt
  - s \lambda \int_{t_0}^{T'}\hspace*{-9pt}\int_{\Omega}
  \nabla  \cdot ({\varphi}
  \nabla \beta) |\psi|^2  \; dx dt
  + 2s \lambda^2 \int_{t_0}^{T'}\hspace*{-9pt}\int_{\Omega}
  {\varphi}
   \vert \nabla \beta \vert^2 |\psi|^2 dx dt,
\end{eqnarray*}
by integration by parts.
With an integration by parts w.r.t. $t$ in the first integral, we
then obtain
\begin{eqnarray*}
  \label{eq:I}
  \frac{1}{2} \int_{\Omega} |\psi(T',.)|^2 \; dx = \mathcal{I}
  - s \lambda^2 \int_{t_0}^{T'}\hspace*{-9pt}\int_{\Omega}
  {\varphi}
  \vert \nabla \beta \vert^2 |\psi|^2 dx
  dt + s \lambda \int_{t_0}^{T'}\hspace*{-9pt}\int_{\Omega}{\varphi}(\Delta \beta )|\psi|^2\; dx dt
\end{eqnarray*}
since $\psi(t_0)=0$ and $\nabla\varphi=\lambda\varphi\nabla\beta$. 
\par
\indent 
Then, we have
\begin{equation}
\label{eq:I2}
  \int_{\Omega}
  (e^{{-2s \eta})(T',x)} |y(T',x)|^2 \; dx
    \leq  2 |\mathcal{I}|
    + C s \lambda(\lambda+1) \int_{t_0}^{T'}\hspace*{-9pt}\int_{\Omega}
    e^{{-2s \eta}(t,x)}{\varphi} |y|^2 dx dt.
\end{equation}
Using \, $\varphi \leq \frac{T^4}{4} \varphi^3$\, the last term in
(\ref{eq:I2}) is overestimated by the left hand side of (\ref{carl-est})
and this last one is absorbed by the l.h.s. of the inequality obtained
in lemma \ref{lemma1}.\\
If we now observe that
$$y(T',x)=\Delta U(T',x)+a U(T',x)+b V(T',x)+ \gamma \widetilde{v}(T',x),$$
we have
$$  | y(T',x) |^2 \geq \frac{1}{2}| \gamma \widetilde{v}(T',x) |^2 - |\Delta U(T',x)+a U(T',x)+b V(T',x) |^2.$$
\begin{flushright}
\rule{.05in}{.05in}
\end{flushright}
The regularity of the solutions of (\ref{syst2}) allows us to write that
for $n \le 3$, $\partial_t \widetilde{v}$ is an element of $L^2(t_0,T, L^{\infty}(\Omega))$.
So, from $|\widetilde{v}(x,T')| \geq r > 0$, we have
$$\exists\,\, k\in L^2(t_0,T), |\partial_t \widetilde{v}(x,t)|\leq k(t)|\widetilde{v}(x,T')|,\,\, \forall x\in \Omega,\,\, t\in (t_0,T).$$
Hence (\ref{eq:lemma2}) can be written
\begin{eqnarray*}
\int_{\Omega} (e^{-2s \eta})(T^{\prime},x) |\gamma|^2|%
\widetilde{v}(x,T^{\prime})|^2 \; dx \\
\leq C s^{-3/2} \lambda^{-2}\int_{t_0}^T \int_{\Omega}  e^{-2s \eta}
 |\gamma|^2 |k(t)|^2|\widetilde{v}(x,T^{\prime})|^2\;dx \;dt \\
+C s^{11/2} \lambda^{6} \int_{t_0}^{T}\hspace*{-9pt} \int_{\omega}
e^{-2s \eta} \varphi^{7} | z |^2\ d x\ d t \\
+C \int_{\Omega} (e^{{-2s \eta})(T^{\prime},x)} (|\Delta
U(T^{\prime},x)|^2 + |U(T^{\prime},x)|^2 + |V(T^{\prime},x) |^2)
\;dx.
  \end{eqnarray*}

Since $k\in L^2(t_0,T)$ implies that $\displaystyle \int_{t_0}^{T}%
\hspace*{-9pt}|k(t)|^2\;dt\leq k_0<+\infty$. 
For $\lambda$ large enough, the term $(1-C
s^{-3/2} \lambda^{-2}k_0)$ can be made positive:
\begin{equation*}
1-C s^{-3/2} \lambda^{-2}k_0\geq C_2>0.
\end{equation*}
Using the fact
that $e^{-2 s \eta(t,x)} \leq e^{-2 s \eta(T^{\prime},x)} \;\;
\forall x \in \Omega,\ \ \forall t \in (t_0,T)$, 
we deduce that
\begin{eqnarray*}
r^2(1-C s^{-3/2} \lambda^{-2}k_0)\int_{\Omega} e^{-2s \eta (T^{\prime},x)} |\gamma|^2 \;dx \\
\leq C s^{11/2} \lambda^{6} \int_{t_0}^{T}\hspace*{-9pt}
\int_{\omega} e^{-2s \eta} \varphi^{7} | z |^2\ d x\ d t \\
+C \int_{\Omega} e^{-2s \eta(T^{\prime},x)} (|\Delta
U(T^{\prime},x)|^2 + |U(T^{\prime},x)|^2 + |V(T^{\prime},x) |^2)
\;dx,
  \end{eqnarray*}
 where we have also used that $|\widetilde{v}(x,T^{\prime})| \geq r > 0$
in $\Omega$.
Then, by virtue of the properties satisfied by $\varphi$ and
$\eta$, we finally obtain
  \begin{eqnarray}
  \label{eq:sta}
  |\gamma|^2_{L^2(\Omega)} \leq \frac{C}{r^2 C_2} s^{11/2}
\lambda^{6} \int_{t_0}^{T}\hspace*{-9pt} \int_{\omega} | z |^2\ d x\ d t \\
+ \frac{C}{r^2 C_2} \int_{\Omega} (|\Delta U(T^{\prime},x)|^2 +
|U(T^{\prime},x)|^2 + |V(T^{\prime},x) |^2) \;dx.\nonumber
  \end{eqnarray}
  With (\ref{eq:sta}), recalling that $U=u-\widetilde{u}$, $V=v-\widetilde{v}$, $y=\partial_t (u-\widetilde{u})$
  and $z=\partial_t (v-\widetilde{v})$, we have thus obtained the following stability result.
  \begin{theorem}
  \label{th:sta}
  Let $\omega$ be a subdomain of an open set $\Omega$ of $\mathbb{R}^n$, let $a,\ b,\ c,\ d\ \in \Lambda(R)$. Furthermore,
we assume that $\widetilde{u}_0$, $\widetilde{v}_0$ in $H^2(\Omega)$ and the
assumption \ref{ass1} is satisfied. Let ($u$, $v$), ($\widetilde{u}$,
$\widetilde{v}$) be solutions to (\ref{syst1})-(\ref{syst2}). Then there
exists a constant C
  $$C=C(\Omega, \omega, c_0, t_0, T, r, R ) > 0$$
  such that
    \begin{eqnarray}
    \label{res:sta}
    |b-\widetilde{b}|^2_{L^2(\Omega)} \leq
    C |\partial_t v - \partial_t \widetilde{v}|^2_{L^2((t_0,T) \times \omega)}
    + C |\Delta u (T', \cdot)- \Delta \widetilde{u}(T', \cdot)|^2_{L^2(\Omega)}\\ 
    + C | u (T', \cdot)-  \widetilde{u}(T', \cdot)|^2_{L^2(\Omega)}
    + C |v (T', \cdot)- \widetilde{v}(T', \cdot)|^2_{L^2(\Omega)}.\nonumber 
    \end{eqnarray}
  \end{theorem}
  \begin{remark}
  If we assume that $u(T', \cdot)=\widetilde{u}(T', \cdot)$ and 
  $v(T', \cdot)=\widetilde{v}(T', \cdot)$ (such an additional assumption
  is sometimes made, e.g. in \cite{IY:98}), then the stability estimate
  becomes
  $$|b-\widetilde{b}|^2_{L^2(\Omega)} \leq
    C |\partial_t v - \partial_t \widetilde{v}|^2_{L^2((t_0,T) \times \omega)}.$$
    \end{remark}
    With Theorem \ref{th:sta} we have the following uniqueness result
    \begin{corollary}
       \label{cor2}
    Under the same assumptions as in theorem \ref{th:sta} and if
    $$(\partial_t v - \partial_t \widetilde{v})(t,x)=0 \;\;\mbox{ in }\;\; (t_0,T) \times \omega,$$
    $$\Delta u (T', x)- \Delta \widetilde{u}(T',x)=0 \;\;\mbox{ in }\;\;\Omega,$$
    $$v (T', x)- \widetilde{v}(T', x)=0 \;\;\mbox{ in }\;\;\Omega,$$
    Then $b=\widetilde{b}$.
    \end{corollary}
%
%
\section{ A uniqueness and stability estimate for the initial conditions}
%
%
In this section, we use the same method as in \cite{YZ:01} 
to state a stability estimate for the initial conditions
$u_0$, $v_0$. 
The idea is to prove logarithmic-convexity inequality. 
The following method has been used to obtain continuous dependence inequalities in
initial value problems.
If $(y,z)$ is solution of (\ref{syst3}), we introduce $(y_1,z_1)$ and $(y_2,z_2)$ that satisfy 
\begin{equation}
\left \{ \begin{array}{lll}
\label{syst5}
 {\partial}_t y_1= \Delta y_1+ay_1+bz_1+\gamma \partial_t \widetilde{v} & \mbox{in} & Q_0,\\
 {\partial}_t z_1= \Delta z_1+cy_1+dz_1 & \mbox{in} & Q_0,\\
    y_1(t,x)=z_1(t,x)=0  & \mbox{on} & \Sigma_0,\\
    y_1(0,x)=0, & \mbox{in} & \Omega,\\
    z_1(0,x)=0 & \mbox{in} & \Omega,
  \end{array}\right.
  \end{equation}
  and
\begin{equation}
\left \{ \begin{array}{lll}
\label{syst6}
 {\partial}_t y_2= \Delta y_2+ay_2+bz_2 \ & \mbox{in}& Q_0,\\
 {\partial}_t z_2= \Delta z_2+cy_2+dz_2 & \mbox{in} & Q_0,\\
    y_2(t,x)=z_2(t,x)=0  & \mbox{on} & \Sigma_0,\\
    y_2(0,x)=\Delta U(0,x)+a U(0,x)+b V(0,x)+ \gamma \widetilde{v}(0,x) & \mbox{in} & \Omega,\\
    z_2(0,x)=\Delta V(0,x)+cU(0,x)+dV(0,x) & \mbox{in} & \Omega. 
  \end{array}\right.
\end{equation} 
 Then, we have
 \begin{equation}
  \label{som}
  y=y_1+y_2 \quad  \mbox{and}\quad  z=z_1+z_2.
  \end{equation}
In a first step, we give an $L^2$ estimate for $(y_1, z_1)$
%
\begin{lemma} 
Let $a,\ b,\ c,\ d,\ |\partial_t \widetilde{v}|_{L^2(t_0,T)}\ \in \Lambda(R)$. Then there exists a constant
$$C=C(t_0, T', R)>0,$$ such that
\begin{equation}
\label{y_1,z_1}
|y_1(t)|^2_{L^2(\Omega)}+|z_1(t)|^2_{L^2(\Omega)}\leq
C|\gamma|^2_{L^2(\Omega)}, \,\,\quad \quad\,\,  t_0\leq t\leq T' .
\end{equation}
\end{lemma}
{\bf Proof}:\\
We multiply the first (resp. the second) equation of (\ref{syst5}) by $y_1$
(resp. by $z_1$). Then, after integrations by parts 
with respect to the space variable, we obtain
$$\begin{array}{l}
\displaystyle
 \frac{1}{2}\partial_t\int_{\Omega}(|y_1|^2+|z_1|^2)\, dx=
-\int_{\Omega}(|\nabla y_1|^2 + |\nabla z_1|^2)\, dx \vspace{0.5cm}\\
\displaystyle\ + \int_{\Omega}ay_1^2\, dx +
\int_{\Omega}dz_1^2\, dx +\int_{\Omega}(b+c)y_1z_1\, dx + \int_{\Omega}\gamma(\partial_t\widetilde{v})y_1\, dx.
\end{array}$$
We use Cauchy-Schwarz and Young inequalities and we integrate
over $(t_0,t)$ for $t_0\leq t\leq T'$, and we obtain
$$|y_1(t)|^2_{L^2(\Omega)}+|z_1(t)|^2_{L^2(\Omega)}\leq
C_2|\gamma|^2_{L^2(\Omega)}+
C_1\int_{t_0}^{t}(|y_1(s)|^2_{L^2(\Omega)}+|z_1(s)|^2_{L^2(\Omega)})\,ds.$$
The result follows by a Gronwall inequality.
\begin{flushright}
\rule{.05in}{.05in}
\end{flushright}
In a second step, we use a logarithmic-convexity inequality for $(y_2, z_2)$
%
%
            
\begin{lemma} 
Let $a,\ b,\ c,\ d\ \in \Lambda(R)$ and $u_0$, $v_0$, $\widetilde{u}_0$,
$\widetilde{v}_0$ in $H^4(\Omega)$. Then there exist constants $M>0$,
$C=C(R)>0$ and $C_1=C_1(t_0, T', R)>0$ such that
\begin{equation}
\label{y_2,z_2}
|y_2(t)|^2_{L^2(\Omega)}+|z_2(t)|^2_{L^2(\Omega)}\leq C_1 M^{1-\mu(t)}
(|y_2(T')|^2_{L^2(\Omega)}+|z_2(T')|^2_{L^2(\Omega)})^{\mu(t)},
\end{equation}
for $t_0\leq t\leq T'$, where 
$\displaystyle{\mu(t)=\frac{(e^{-Ct_0}-e^{-Ct})}{(e^{-Ct_0}-e^{-CT'})}}$.
\end{lemma}
{\bf Proof}:\\
The proof of this lemma is just an application of Theorem 3.1.3 in
\cite{I:98}.
In fact, system (\ref{syst6}) can be written in the following form
 \begin{eqnarray*}
 {\partial}_t W+ AW=BW,\ & \mbox{in $Q_0$},\\
  W(t,x)=0  & \mbox{on $\Sigma_0$},\\
   W(0,x)=W_0(x) & \mbox{in $\Omega$},\\
  \end{eqnarray*}
  where 
 $$W=
 \left (\begin{array}{c}
   y_2\\
   z_2
 \end{array}\right)
 ,\quad
  A=
 \left ( \begin{array}{cc}
   -\Delta & 0\\
   0 & -\Delta \\
 \end{array}\right)
 ,\quad
 B=
 \left ( \begin{array}{cc}
   a & b\\
   c & d
 \end{array}\right).
 $$
 The operator $A$ is symetric and the solution $W$ satisfies
 $$||\partial_t{W}+AW||_{L^2(\Omega)}\leq \alpha ||W||_{L^2(\Omega)},$$\\
 where $\alpha=||B||_{L^\infty(\Omega)}<+\infty$ since $a$, $b$, $c$ and $d$ 
are in $ \Lambda(R)$.
If we assume that $u_0$, $v_0$, $\widetilde{u}_0$,
$\widetilde{v}_0$ are in $H^4(\Omega)$, 
the hypothesis of Theorem 3.1.3 in \cite{I:98} are satisfied, thus we have
 $$ ||W(t)||_{L^2(\Omega)}\leq C_1 ||W(t_0)||^{1-\mu(t)}_{L^2(\Omega)}
||W(T')||^{\mu(t)}_{L^2(\Omega)},$$
 with $\displaystyle{\mu(t)=\frac{(e^{-Ct_0}-e^{-Ct})}{(e^{-Ct_0}-e^{-CT'})}}$.\\
 \noindent
Since $W \in {\cal C}(t_0,T, L^2(\Omega))$, we have
 $||W(t_0)||_{L^2(\Omega)}\leq M^{\frac{1}{2}}$, and the result follows.
\begin{flushright}
\rule{.05in}{.05in}
\end{flushright}
 The two previous lemmas allow us to prove the following 
%
\begin{theorem}
 Let $\omega$ be a subdomain of an open set $\Omega$ of $\mathbb{R}^n$, let $a,\ b,\ c,\ d\ \in \Lambda(R)$.
 Furthermore,
we assume that $u_0$, $v_0$, $\widetilde{u}_0$, $\widetilde{v}_0$ in
$H^4(\Omega)$ and Assumption \ref{ass1} is satisfied. 
Let
($u$, $v$), ($\widetilde{u}$, $\widetilde{v}$) be solutions to
(\ref{syst1})-(\ref{syst2}). We set $$ \,\, E=|\partial_t v - \partial_t
\widetilde{v}|^2_{L^2((t_0,T) \times
\omega)}+|u(T', \cdot)-\widetilde{u}(T', \cdot)|^2_{H^2(\Omega)}+|v(T', 
\cdot)-\widetilde{v}(T', \cdot)|^2_{H^2(\Omega)}.$$
Then there exists a constant $C=C(\Omega, \omega, c_0, t_0, T, r, R ) > 0$
  such that
$$
|u_0-\widetilde{u}_0|^2_{L^2(\Omega)}+|v_0-\widetilde{v}_0|^2_{L^2(\Omega)}\leq 
\frac{C}{|\log E|},\;\; \mbox{ for } \;\; 0<E<1
$$
\end{theorem}
{\bf Proof}:\\
Since $\widetilde{v} \in H^1(t_0,T, H^2(\Omega))$, we have $\widetilde{v} \in L^{\infty}(\Omega)$.
In view of (\ref{som}), inequalities (\ref{y_1,z_1}), (\ref{y_2,z_2}) imply
\begin{eqnarray*}
|y(t, \cdot)|^2_{L^2(\Omega)}\leq 2(|y_1(t, \cdot)|^2_{L^2(\Omega)} 
+ |y_2(t, \cdot)|^2_{L^2(\Omega)})\\
\leq C_1 \left [|\gamma|^2_{L^2(\Omega)} + M^{1-\mu(t)}(|y_2(T', \cdot)|^2_{L^2(\Omega)}
+|z_2(T', \cdot)|^2_{L^2(\Omega)})^{\mu(t)} \right] .\nonumber 
\end{eqnarray*}
Now, with (\ref{som}), we write $y_2=y-y_1$ and $z_2=z-z_1$.
Inequality (\ref{y_1,z_1}) gives us an estimation of 
$|y_1(T', \cdot)|^2_{L^2(\Omega)}$ and
$|z_1(T', \cdot)|^2_{L^2(\Omega)}$ in terms of
$ |\gamma|^2_{L^2(\Omega)}$. Then the definition of $y$ and $z$ in
(\ref{syst3}) gives us an estimation of 
$|y(T', \cdot)|^2_{L^2(\Omega)}$ and
$|z(T', \cdot)|^2_{L^2(\Omega)}$ in terms of 
$|U(T', \cdot)|^2_{H^2(\Omega)}$ and $|V(T',\cdot)|^2_{H^2(\Omega)}$.
Finally, we obtain
$$
|y(t, \cdot)|^2_{L^2(\Omega)}
\leq C_1 \left[|\gamma|^2_{L^2(\Omega)} + M_2(|\gamma|^2_{L^2(\Omega)} + 
|U(T', \cdot)|^2_{H^2(\Omega)}+|V(T',\cdot)|^2_{H^2(\Omega)})^{\mu(t)} \right],
$$
(a similar estimate is obtained for $|z(t, \cdot)|^2$).
Following \cite{YZ:01} with a time translation, if we use (\ref{res:sta}), the last estimate yields
\begin{eqnarray*}
|u_0-\widetilde{u}_0|^2_{L^2(\Omega)}+|v_0-\widetilde{v}_0|^2_{L^2(\Omega)}=|U(t_0,
\cdot)|^2_{L^2(\Omega)} + |V(t_0, \cdot)|^2_{L^2(\Omega)}\\
=|-\int_{t_0}^{T'}y(s, \cdot)\ ds + U(T',
\cdot)|^2_{L^2(\Omega)} + |-\int_{t_0}^{T'}z(s, \cdot)\ ds + V(T',
\cdot)|^2_{L^2(\Omega)}\nonumber \\
\leq M_3\int_{t_0}^{T'} E^{\mu(s)}\ ds + C_2E
\leq \displaystyle{C_3\frac{E-1}{\log E}}+C_4E
\leq \displaystyle{\frac{C}{|\log E|}}.\nonumber 
\end{eqnarray*}\begin{flushright}
\rule{.05in}{.05in}
\end{flushright}
\ack
The authors wish to thank A. Benabdallah, Y. Dermenjian, M. Henry, J. Le Rousseau 
and O. Poisson for numerous usefull discussions on the subject of this paper. 

\section*{References}

\end{document}